\documentclass[11pt]{amsart}

\usepackage[textheight=815 pt, textwidth=360 pt]{geometry}
\usepackage{amssymb}

\usepackage{hyperref}

\newtheorem{thm}{Theorem}[section]
\newtheorem{prop}[thm]{Proposition}
\newtheorem{lem}[thm]{Lemma}
\newtheorem{cor}[thm]{Corollary}
\newtheorem{conj}[thm]{Conjecture} 

\theoremstyle{definition}
\newtheorem{definition}[thm]{Definition}

\theoremstyle{remark}

\numberwithin{equation}{section}

\begin{document}

\title{Charles Bouton and the Navier-Stokes Global Regularity Conjecture}
\author{J. Polihronov}
\address{KLS Inc., 1001 Fanshawe College Blvd. Room A1029, London, Ontario N5Y 5R6, Canada}
\email{j.polihronov@gmail.com}
\urladdr{https://www.researchgate.net/profile/J-Polihronov}


\begin{abstract}
This article examines the Bouton-Lie group invariants of the Navier-Stokes equation (NSE) for incompressible fluids. Bouton's theory is applied to the general scaling transformation admitted by the NSE and is used to derive all self-similar solutions. In light of these, the criticality of the standard NSE system is examined and criticality criteria are derived. The theorem of Beale-Kato-Majda is used to rule out blow-up for a subset of Bouton's self-similar solutions. For a subset of Leray's self-similar solutions, the cavitation number of the fluid is found to be a scale-invariant, conserved quantity. By extending the analysis of Bouton to higher-dimensioned manifolds, additional conserved quantities are found, which could further elucidate the physics of fluid turbulence.
\end{abstract}

\maketitle

\section{Introduction}
Global regularity of the incompressible Navier-Stokes equation (NSE) is subject of intense research. The present article adds to this effort a study of the symmetry properties of the NSE by following the example of a similar work on the general linear differential equation (GLDE, given below) by Charles L. Bouton \cite{Bou99}. We present the work of Bouton on the GLDE and then apply his analysis to the NSE. Following the reasoning of Bouton, we examine the invariants of the NSE with the purpose to draw conclusions about the solutions and their properties. 

The NSE is given by
\begin{eqnarray}
	\rho \left( \frac{\partial \vec{u}}{\partial t} + (\vec{u} \cdot \nabla)\vec{u} \right) &=& \mu \Delta \vec{u} - \nabla p \nonumber\\
	\nabla \cdot \vec{u} &=& 0 \nonumber\\
	\vec{u}(x,y,z,0)&=&\vec{u}_0(x,y,z)
	\label{NSE}
\end{eqnarray}
We suppose initial conditions $\vec{u}_0$ in the Schwartz class \cite{Cla06}. Here $\vec{u}(x,y,z,t)=(u,v,w)$ is the velocity, $\mu=\rho \nu$ is the dynamic viscosity, $\nu$ is the kinematic viscosity, $\rho=\mbox{const}$ is the fluid density and $p(x,y,z,t)$ is the pressure.
\begin{conj}
The NSE global regularity conjecture \cite{Cla06} is the assumption that given the above initial conditions, there always exists a solution $\vec{u}(x,y,z,t)$, $p(x,y,z,t)$ which remains smooth of class $C^{\infty}$ throughout ${\rm I\!R}^3$ and  for all $t \ge 0$:
\begin{eqnarray*}
	\vec{u}_i(x,y,z,t) \in C^{\infty}({\rm I\!R}^3 \times [0,\infty)) \nonumber \\
	p(x,y,z,t) \in C^{\infty}({\rm I\!R}^3 \times [0,\infty)) \nonumber \\
	i = 1, 2, 3; \mbox{   } u_1 = u, u_2 = v, u_3 = w,
\end{eqnarray*}
as well as has bounded energy
\begin{equation*}
	\int _{\rm I\!R^3}  |\vec{u}(x,y,z,t)|^2 dxdydz < C, \forall t \ge 0.
\end{equation*}
	\label{Conjecture}
\end{conj}
The NSE assumes continuous fluid and thus its velocity exists at infinitely small points. This presents difficulty, because despite that the total kinetic energy is conserved, the scaling properties of the NSE prevent it from being able to control the solution at very small scales. This opens the possibility for scaling induced blow-up and the solution may no longer be smooth after some finite time. If it can be demonstrated that the solution can remain smooth this will further justify the choice of Navier and Stokes to model a fluid as a continuous medium and will show that the matematics of their equation always has physically meaningful solutions.

We set about to work towards this goal by studying the symmetry properties of the NSE. Parallel work was completed on the GLDE by Charles L. Bouton as early as 1899 \cite{Bou99}. Bouton examines in detail Lie's invariants of the GLDE. His is one of the earliest works where a differential equation is thoroughly studied with Lie's theory of continuous groups. Due to its comprehensive analysis, Bouton's article can be considered foundational in invariant theory \cite{Soh16}. Historically, the theory of point transformations and Lie groups has been a well-known tool for the study of differential equations \cite{Coh11}. These techniques are now widely used in the analysis of PDEs and are found in the textbooks \cite{Olv93}. To the author's best knowledge, they have not been applied previously to study existence and smoothness of the solutions of differential equations and are used to this end here for the first time.

This article is organized as follows: in the first section, we present the work of Bouton on the GLDE; then, we apply his analysis to the NSE. We use the most general scaling transformation, as published by Ercan-Kavvas \cite{Erc15}, wherein the viscosity $\nu$ is allowed to scale-transform. Bouton's self-similar solutions are derived by applying his first theorem. Next, we re-examine NSE criticality based on the derived solutions. We show that, under the general scaling transformation, the energy of the standard NSE system is no longer supercritical and can be critical or subcritical. Criticality criteria are derived and shown to be governed by a 5/2 law. Using the Beale-Kato-Majda criterion \cite{BKM}, we then argue that scaling induced blow-up will not take place when the NSE problem is energy-subcritical. Bounded energy is shown to follow from the same line of reasoning. 

As a next step, we examine the scaling transformation where the viscosity $\nu$ is not allowed to transform but is kept fixed at all length scales. Here, Bouton's theory leads to the self-similar solutions of Leray \cite{Erc15}, \cite{Brad18}. The aforementioned solutions of Bouton in the case of this transformation are non self-similar polynomials. We also show that for a subset of Leray's solutions we can derive a scale-invariant, conserved quantity, namely Euler's or, cavitation number of the fluid. 

In the last sections, we employ Cartan's calculus in order to show that Bouton's first theorem is but the first of nine steps in the differential geometry analysis of the NSE. In each stage, the Lie derivatives of differential forms are set to zero along the flows of the infinitesimal operators admitted by the NSE and are solved as a system, which yields their common invariants. By virtue of Bouton's theorems, selected invariants are shown to be conserved quantities. Seven such conserved quantities are found; they are differential forms on manifolds, associated with the flow of liquids. Their role as conditions on the pressure and velocity of the fluid in turbulent flows is discussed.

\section{Bouton and the GLDE}
Charles Leonard Bouton is an American mathematician who has conducted a detailed study of the invariant properties of the GLDE \cite{Bou99}. We begin by quoting relevant paragraphs from Bouton's article on the GLDE \cite{Bou99}; we start at \S13. The GLDE is given by 
\begin{equation*}
	y^{(n)}+na_1 y^{(n-1)}+\frac{n(n-1)}{2}a_2y^{(n-2)}...+a_ny=0,
\end{equation*}
where $a_s$ and $y$ are functions of $x$ alone and $a_0 \equiv 1$. Also given is the point transformation
\begin{eqnarray}
	x_1&=&\chi(x) \nonumber \\
	y_1&=&y\psi(x).
	\label{GLDEgroup}
\end{eqnarray}
known as {\it finite} transformation. Following Bouton's notation, $x_1$ and $y_1$ are the transformed variables; $\chi$ and $\psi$ are arbitrary analytical functions.

\S13: A differential equation is said to admit a point transformation, if this transformation leaves the equation unchanged, or, transforms it into an equation of the same form. The transformation (\ref{GLDEgroup}) is the most general point transformation the GLDE is known to admit. Next, Bouton introduces the following definitions:

\begin{definition}
An {\it absolute invariant/covariant} $I$ of the GLDE is a function of the coefficients $a_s$ and their derivatives. It remains unchanged upon transformation (\ref{GLDEgroup}), $I_1 = I$. Invariants do not contain the dependent variable $y$ or its derivatives; covariants do (as per \S26). Bouton probably introduces the term {\it covariant}, since the dependent variable $y$ co-varies with $a_s$. $I$ may contain the independent variable $x$.
\end{definition}

\begin{definition}
	A {\it relative invariant/covariant} $R$ (as per \S38) of the GLDE is a function of the coefficients $a_s$ and their derivatives. Upon finite transformation (\ref{GLDEgroup}) $R$ transforms according to the rule $R_1=f R$; $f$ is known as the factor of finite transformation. Invariants do not contain the dependent variable $y$ or its derivatives; covariants do (as per \S26). $R$ may contain the independent variable $x$.
\end{definition}

\begin{definition}
Define Isobaric weights (\S15): $[y^{(\mu)}]^\nu$ has weight $\mu\nu$; $[a_i^{(j)}]^l$ has weight $(i+j)l$. 
\end{definition}

Bouton begins deriving the invariants of the GLDE as early as chapter 2 and continues in chapters 4 and 5. He uses the method of infinitesimal transformations by solving the equation 
\begin{equation}
	XI = 0,
	\label{XI}
\end{equation}
where $X$ is the symbol/operator of infinitesimal transformation, corresponding to the finite transformation (\ref{GLDEgroup}). Equation (\ref{XI}) represents the necessary and sufficient condition for $I$ to be an absolute invariant/covariant of the GLDE. In (\ref{XI}), $X$ is not only an operator of the geometric transformation (\ref{GLDEgroup}), but is also extended \cite{Bou98} to include terms from the GLDE in itself: $X=X(x,y,\chi,\psi,a_s,a_i^{(j)},y^{(\mu)})$. 

\S15: Write the arguments of the absolute covariant $I$ (they are the same for the relative covariant $R$): 
\[
	I=I(x,y,y',y'',...,y^{(\mu)},a_1,a_2,...,a_i^{(j)},...).
\]

\S14,\S15: Consider the two scaling transformations
\begin{eqnarray}
	x_1 = x  & \mbox{   } & x_1 = Cx \nonumber\\
	y_1 = Cy & \mbox{   } & y_1 = y \nonumber\\
	y_1^{(\mu)}=Cy^{(\mu)} & \mbox{   } & y_1^{(\mu)}=C^{-(\mu)} y^{(\mu)}  \nonumber\\
	\psi=C & \mbox{   } &  \psi=1  \nonumber\\
	\chi=x & \mbox{   } &  \chi=Cx \nonumber\\ 
	& C=\mbox{const}. 
	\label{GLDEscaling}
\end{eqnarray}
Since these are included in (\ref{GLDEgroup}) they influence the properties of all relative invariants and covariants $R$ of the GLDE. When $R$ transforms, it must be according to the law $R_1=fR$. Therefore, $R$ must be an integral rational function (IRF), isobaric of weight $W$ and homogeneous of degree $\lambda$.

\S14, \S15: The scaling transformations (\ref{GLDEscaling}) also determine the properties of all absolute invariants and covariants $I$ of the GLDE. In order to have $I$ transform as $I_1=I$, it follows that $I$ must be an IRF of weight zero and also homogeneous of degree zero. Bouton thus arrives at the following theorem:

\begin{thm}
(of Bouton): Any absolute invariant/covariant of the GLDE for the group of transformations (\ref{GLDEgroup}) must be homogeneous in the $y^{(\mu)}$'s of degree zero and isobaric in the $y^{(\mu)}$'s and $a_i^{(j)}$'s of weight zero. Any relative invariant/covariant $R$ must be homogeneous in the $y^{(\mu)}$'s of degree $\lambda$ and isobaric in the $y^{(\mu)}$'s and $a_i^{(j)}$'s of weight $W$. We write this as $R^{(\lambda,W)}$.
\end{thm}

\begin{cor}
The transformation of a function is determined by the transformable entities within its functional form. This is seen in the reasoning of Bouton by the way he calculates isobaric weights and powers of homogeneity. Under finite transformation, the function $I(x,y,y',y'',...,y^{(\mu)},a_1,a_2,...,a_i^{(j)},...)$ will transform according to how its individual arguments would transform within the functional form $I$. The same is true for the relative invariants/covariants $R$. 
	\label{Cor1-Bouton}
\end{cor}

\begin{cor}
In the most general case, the factor of the relative invariants/covariants of the GLDE is $f=\psi^\lambda/\chi'^W$ (see eq.(\ref{GLDEgroup})).
\label{Cor2-Bouton}
\end{cor}

\begin{cor}
Due to the requirements for isobarity and homogeneity, IRFs are the only functions arising in this study of invariants. In the particular cases when the denominator is constant, the IRF is a homogeneous, isobaric polynomial (e.g. $y^{(12)}+(y^{\prime\prime\prime})^4+(y^{\prime\prime})^6$). 
	\label{Cor3-Bouton}
\end{cor}

\section{Bouton and the NSE}
The finite transformations and their corresponding infinitesimal group operators which leave the NSE invariant have been studied and published extensively in the literature. They were derived in 3-dimensional space by Lloyd \cite{Llo81}: the NSE admits time translations, scale changes and rotations (details provided in later chapters), the finite scaling transformation being
\begin{eqnarray}	
	(x',y',z')&=&k(x,y,z) \nonumber \\
	t'&=&k^2 t \nonumber \\
	(u',v',w')&=&(1/k)(u,v,w) \nonumber \\
	p'&=&(1/k^2)p,  \nonumber \\
	0<k<\infty,& &k=\mbox{const}, 
\label{scaling}
\end{eqnarray}
where primed quantities are transformed quantities; we no longer use the index ``1" to denote transformation. We omit the admitted groups that are not Lie groups. The NSE does not admit any other groups \cite{Llo81}. According to recent work by Ercan and Kavvas the general scaling transformation admitted by the Navier-Stokes equation is \cite{Erc15}:

\begin{eqnarray}	
	(x',y',z')&=&k^{\alpha_x}(x,y,z) \nonumber \\ 
	t'&=&k^{\alpha_t} t \nonumber \\
	\rho^{\prime} &=& k^{\alpha_\rho} \rho  \nonumber \\
	(u',v',w')&=&k^{\alpha_x-\alpha_t}(u,v,w) \nonumber\\
	p'&=&k^{\alpha_\rho+2\alpha_x-2\alpha_t}p  \nonumber \\
	\nu^{\prime} &=& k^{2\alpha_x-\alpha_t} \nu  \nonumber \\
	g_i' &=& k^{\alpha_x-2\alpha_t} g_i, i=1,2,3  \nonumber \\
	0< k< \infty,& & k=\mbox{const} 
      	\label{EK-scaling}
\end{eqnarray}
and if $\alpha_t=2, \alpha_x=1$ with gravitational acceleration $g_i=0$, we arrive at the scaling transformation (\ref{scaling}). In the general case (\ref{EK-scaling}), all three spatial variables transform with the same scaling exponent $\alpha_x$ and all three velocity components transform with the same scaling exponent $\alpha_x-\alpha_t$; $\alpha_x$ and $\alpha_t$ are arbitrary real numbers. 

In this work, we set $g_i=0$. Also, we set ``$p$" to stand for $p/\rho$ and thus the finite transformation for $p$ in (\ref{EK-scaling}) will become $p'=k^{2\alpha_x-2\alpha_t}p$. If the NSE is dimensionalized, $\nu$ must scale-transform since it is a dimensional constant ($\nu$ has dimension) and may be thought of $\nu=x_0^2/t_0$. In the nondimensionalized NSE the coefficient $\nu$ no longer denotes viscosity, but the reciprocal value of the Reynolds number $\mbox{Re}^{-1}$; if Re is not allowed to scale, (\ref{EK-scaling}) still applies as long as $\alpha_t=2\alpha_x$. However, as we will see later, Leray's self-similar solutions are the solutions of the nondimensionalized NSE and their self-similarity is due to scaling invariance; obviously all variables in the nondimensionalized NSE are unitless and yet allowed to scale. We will therefore examine in more detail both cases for the viscosity $\nu$: (1) $\nu$ is allowed to scale; (2) $\nu$ is kept fixed at all length scales. Analogously with Bouton's definitions of the preceding section, we introduce

\begin{definition}
	An {\it absolute invariant/covariant} $I$ of the NSE is a function of the independent variables $x,y,z,t$. $I$ is such that upon transformations (\ref{scaling})-(\ref{EK-scaling}) it remains unchanged, $I' = I$. Invariants do not contain the dependent variables $u,v,w,p$ or their derivatives; covariants do.
\end{definition}

\begin{definition}
	A {\it relative invariant/covariant} $R$ of the NSE is a function of the independent variables $x,y,z,t$. $R$ transforms according to the rule $R'=f R$  upon transformation (\ref{scaling})-(\ref{EK-scaling}); $f$ is the transformation factor. Invariants do not contain the dependent variables $u,v,w,p$ or their derivatives; covariants do.
\label{rel-cov}
\end{definition}

\begin{definition}
Define weights as follows: 
\begin{eqnarray*}
	W(x^a)=W(y^a)=W(z^a)&=&a\alpha_x \nonumber \\
	W(t^b)&=&b\alpha_t \nonumber\\
	W(u^a) = W(v^a)=W(w^a)&=&(\alpha_x-\alpha_t)a \nonumber \\ 
	W(p^b)&=&2(\alpha_x-\alpha_t)b \nonumber\\
	W(\nu^c)&=&(2\alpha_x-\alpha_t)c.
\end{eqnarray*}
When variables are multiplied, their weights are added. 
\end{definition}

Consider the finite scaling transformation: the relative invariants and covariants $R$ of the NSE must transform according to the law $R'=fR$. From Corollary \ref{Cor2-Bouton} of Bouton, Definition \ref{rel-cov} and from eqs. (\ref{scaling})-(\ref{EK-scaling}) it follows that here $f=k^q$. If $\alpha_t=2,\alpha_x=1$ the relative covariants $R=u$, $R=v$ and $R=w$ yield $f=k^{-1}$, while the relative covariant $R=p$ yields $f=k^{-2}$. Corollary \ref{Cor1-Bouton} of Bouton states that under finite transformation the dependent variables (the jet-space variables, or relative covariants) $u$, $v$, $w$ and $p$ transform the same way as the relative invariant functions $u(x,y,z,t)$, $v(x,y,z,t)$, $w(x,y,z,t)$ and $p(x,y,z,t)$. To ensure that their factor $f$ is correct upon scaling transformation, they must be IRFs; they must be isobaric of their respective weight in the $x,y,z$ and also in $t$. Note that as in Corollary \ref{Cor3-Bouton} of Bouton, IRFs are the only functions arising in the study of the invariants of the NSE. Whenever the denominator of the IRF is constant, the invariant is an isobaric polynomial, e.g. $R^{(2,-2)}=u^2+v^2+w^2$. All presented reasoning in this section therefore sums up in the following theorem:

\begin{thm}
(of Bouton): To meet all isobarity requirements of their scaling transformation, the relative invariant functions $u(x,y,z,t)$, $v(x,y,z,t)$ and $w(x,y,z,t)$ must be isobaric in $x,y,z$ and $t$ of weight $W=\alpha_x-\alpha_t$; while the relative invariant function $p(x,y,z,t)$  must be isobaric in $x,y,z$ and $t$ of weight $W=2(\alpha_x-\alpha_t)$. The Navier-Stokes system of PDEs becomes: 
\begin{eqnarray}
	\frac{\partial \vec{u}}{\partial t} + (\vec{u} \cdot \nabla)\vec{u} &=& \nu \Delta \vec{u} - \nabla p \nonumber\\
	\nabla \cdot \vec{u} &=& 0 \nonumber\\
	(\vec{r}\cdot\nabla)\vec{u}+t\frac{\alpha_t}{\alpha_x}\frac{\partial \vec{u}}{\partial t}&=&\frac{\alpha_x-\alpha_t}{\alpha_x}\vec{u} \nonumber\\
	(\vec{r}\cdot\nabla)p+t\frac{\alpha_t}{\alpha_x}\frac{\partial p}{\partial t}&=&\frac{2(\alpha_x-\alpha_t)}{\alpha_x} p.
	\label{NSE1}	
\end{eqnarray}
\label{T2-Bouton}
\end{thm}
Here, $\vec{r}=(x,y,z)$ is the position vector. The added equations represent Bouton's isobarity conditions. 

\begin{cor}
The transformation of a function is determined by the transformable entities within its functional form. This is seen in the reasoning of Bouton by the way he calculates isobaric weights and powers of homogeneity. Under finite transformation, the function $I(x,y,z,t)$ will transform according to how its individual arguments would transform within the functional form $I$. The same is true for the relative invariants $\vec{u}(x,y,z,t),p(x,y,z,t)$. 
\end{cor}

\begin{cor}
	Under the general scaling transformation (\ref{EK-scaling}), the factor of the relative invariants/covariants of the NSE is $f=k^{\alpha_x-\alpha_t}$ for the velocity and $f=k^{2(\alpha_x-\alpha_t)}$ for the pressure.
\end{cor}

\begin{cor}
	Due to the requirements for isobarity (and also homogeneity, see later chapters), IRFs are the only functions arising in this study of the invariants of the NSE. In the particular cases when the denominator is constant, the IRF is an isobaric polynomial, for ex. $(t^3/t_0^4)(x+z+(y^2/t_0^4))$; $t_0=$ const. 
\end{cor}

\section{When Viscosity is allowed to scale}
To find Bouton's self-similar solutions we apply Theorem \ref{T2-Bouton}; with the symbolic package Maple \cite{Map19}, the added equations of Bouton can be integrated to yield:
\begin{eqnarray}
	\vec{u}& =& t^{\frac{\alpha_x-\alpha_t}{\alpha_t}} \mathbf{F} \left( \frac{x}{t^{\frac{\alpha_x}{\alpha_t}}},\frac{y}{t^{\frac{\alpha_x}{\alpha_t}}} ,\frac{z}{t^{\frac{\alpha_x}{\alpha_t}}}  \right)  \nonumber \\
	p &=& t^{\frac{2(\alpha_x-\alpha_t)}{\alpha_t}} F \left( \frac{x}{t^{\frac{\alpha_x}{\alpha_t}}},\frac{y}{t^{\frac{\alpha_x}{\alpha_t}}} ,\frac{z}{t^{\frac{\alpha_x}{\alpha_t}}}  \right).
	\label{gen-selfsimilarsol}
\end{eqnarray}
These expressions are the self-similar solutions of the PDE system (\ref{NSE1}); they are isobaric IRF's due to the scaling invariance of the NSE. The general form (\ref{gen-selfsimilarsol}) is a necessary condition that all self-similar solutions must meet. 

As the NSE models a fluid continuum, solutions will exist at increasingly small scales. We need to address the possibility of scaling-induced blow-up, where mechanical energy may be stored within infinitely small volume while maintaining energy conservation. The NSE problem may be supercritical, critical or subcritical depending on how its conserved quantities scale when the solution is subjected to a scaling transformation \cite{Tao07}. As an example, take the energy in the case $\alpha_t=2, \alpha_x=1$ for the NSE in $D$ dimensions: the rescaling that preserves the equation is given in (\ref{scaling})
\begin{equation*}
   \vec{u}(x,y,z,t) \rightarrow \frac{1}{k}\vec{u}\left( \frac{x}{k},\frac{y}{k},\frac{z}{k},\frac{t}{k^2} \right),
\end{equation*}
while the energy is an integral of $|\vec{u}|^2$ and thus the rescaling that preserves the energy is given by
\begin{equation*}
   \vec{u}(x,y,z,t) \rightarrow \frac{1}{k^{D/2}} \vec{u}\left(\frac{x}{k},\frac{y}{k},\frac{z}{k},\frac{t}{k^2}\right),
\end{equation*}
thereby showing that for the NSE in $D=3$ dimensions energy is being rescaled more severely; this amplifies the rescaling at small scales $k\ll1$ \cite{Tao07}. Thus, the case $\alpha_t=2$, $\alpha_x=1$ is energy-supercritical in 3-D. Similar procedure for establishing criticality is found in \cite{Tao08}; suppose $X$ is a conserved quantity; then, if $\lambda$ is a scaling coefficient,

(a) the PDE is $X$-critical, if $X^\prime=X$ $\forall \lambda$ (scale-invariant);

(b) the PDE is $X$-subcritical, if $X^\prime$ diminishes as $\lambda$ increases ($\lambda > 1$);

(c) the PDE is $X$-supercritical, if $X^\prime$ increases as $\lambda$ increases ($\lambda > 1$).

To ensure meaningful physics of the fluid at $t=0$ and throughout $\rm I\!R^3$, we formulate the following lemma for Bouton's self-similar solutions (\ref{gen-selfsimilarsol}):

\begin{lem}
Bouton's self-similar solutions (\ref{gen-selfsimilarsol}) can be nonzero and smooth at $t=0$ if and only if they are in the form
\begin{eqnarray}
	\vec{u} &=&\mathbf{C} t^{\frac{\alpha_x-\alpha_t}{\alpha_t}} +  x^{\frac{\alpha_x-\alpha_t}{\alpha_x}} \mathbf{P} \left( \frac{y}{x},\frac{z}{x}  \right)  \nonumber \\
	p &=& C t^{\frac{2(\alpha_x-\alpha_t)}{\alpha_t}}  +   x^{\frac{2(\alpha_x-\alpha_t)}{\alpha_x}} P \left( \frac{y}{x},\frac{z}{x}  \right),
	\label{gen-selfsimilarsol1}
\end{eqnarray}
where $\mathbf{C}$ and $C$ are constants and the remaining terms are isobaric polynomials.
\label{Lemma}
\end{lem}

\begin{proof}
	Bouton's self-similar solutions are in the form (\ref{gen-selfsimilarsol}). Based on the symmetry properties of the NSE, we will assume that isobaric solutions (i.e. those that scale-transform according to (\ref{EK-scaling})) always exist. To this assumption, we add the Conjecture (\ref{Conjecture}), namely that our chosen solution is smooth at $t=0$ in order to guarantee meaningful value of the initial energy. Smoothness of (\ref{gen-selfsimilarsol}) at $t=0$ requires $(\alpha_x-\alpha_t)/\alpha_t >0$ and  $(\alpha_x-\alpha_t)/\alpha_x >0$, that is, smoothness of the arbitrary functions $\mathbf{F}$ and $F$ everywhere in $\rm I\!R^3$ and for $t=0$. If $t$ is a factor in (\ref{gen-selfsimilarsol}) however, even though smooth, such solution is trivial as it is identically zero when $t=0$; and since there is no forcing term in the NSE,  eq.(\ref{NSE}), the fluid system would be left with a total energy of zero, because its sole energy source is the initial velocity of the system. For this reason, we add the arbitrary isobaric functions $\mathbf{P}$ and $P$ (they are still in the isobaric form (\ref{gen-selfsimilarsol})),  which must themselves be smooth everywhere in $\rm I\!R^3$:

\begin{eqnarray*}
	\vec{u} & = & t^{\frac{\alpha_x-\alpha_t}{\alpha_t}} \mathbf{F} \left( \frac{x}{t^{\frac{\alpha_x}{\alpha_t}}},\frac{y}{t^{\frac{\alpha_x}{\alpha_t}}} ,\frac{z}{t^{\frac{\alpha_x}{\alpha_t}}}  \right) +  x^{\frac{\alpha_x-\alpha_t}{\alpha_x}} \mathbf{P} \left( \frac{y}{x},\frac{z}{x}  \right)  \nonumber \\
	p &=& t^{\frac{2(\alpha_x-\alpha_t)}{\alpha_t}} F \left( \frac{x}{t^{\frac{\alpha_x}{\alpha_t}}},\frac{y}{t^{\frac{\alpha_x}{\alpha_t}}} ,\frac{z}{t^{\frac{\alpha_x}{\alpha_t}}} \right) +   x^{\frac{2(\alpha_x-\alpha_t)}{\alpha_x}} P \left( \frac{y}{x},\frac{z}{x}  \right).
\end{eqnarray*}

Regarding smoothness at $t=0$: $\mathbf{F}$ and $F$ must be smooth in $({\rm I\!R}^3 \times 0)$; however, they are IRFs and as such must be either isobaric polynomials, or the ratio of such polynomials (note that such polynomials cannot contain constant terms since they do not have isobaric weight). Suppose they are a ratio of isobaric polynomials; then their denominators must be nonzero at $t=0$. However, at $t=0$ their denominators must be isobaric polynomials of only $x,y,z$ which are always zero at least at the origin. This is a contradiction with the requirement of smoothness. Therefore, $\mathbf{F}$ and $F$ must be isobaric polynomials. However, it is impossible to build a smooth polynomial by using the arguments $x_i/t^{(\alpha_x/\alpha_t)}$; since multiplying or adding them always leads to discontinuities in either $\rm I\!R^3$ (esp. at the origin) or at $t=0$. Similarly, $\mathbf{F}$ and $F$ cannot be closed-form functions of the given arguments, since the arguments are discontinuous in $({\rm I\!R}^3 \times 0)$. Therefore, the only remaining option is that they are constants.

	Regarding the smoothness of $x^{\frac{\alpha_x-\alpha_t}{\alpha_x}} \mathbf{P}$ and $x^{\frac{2(\alpha_x-\alpha_t)}{\alpha_x}} P$: in order for the remaining terms to be smooth in $\rm I\!R^3$ (as per the requirement of smoothness at $t=0$), they must be isobaric polynomials. They cannot be IRFs, since their denominators will be polynomials with roots at the origin, etc. They cannot be closed-form functions, since their Taylor series will be either isobaric but discontinuous at least at the origin, or not isobaric. Finally, it is evident that $(\alpha_x-\alpha_t)/\alpha_x >0$ is mandatory to ensure smoothness.
\end{proof}

Therefore, according to Lemma \ref{Lemma}, if the self-similar solution is to be smooth and nonzero at $t=0$, the following condition must be true:
\[
	\frac{\alpha_x-\alpha_t}{\alpha_t} > 0, 
\]
that is, the power index of $t$ must be positive; and it is the stronger requirement in comparison to the requirement of the power index of $x$ being positive. Lemma \ref{Lemma} leads to the following corollaries

\begin{cor}
	The power index of both $t$ and of $x$ in (\ref{gen-selfsimilarsol1}) is positive if 
\begin{align*}
	& \alpha_x > \alpha_t \mbox{ and } \alpha_t > 0, \mbox{ thus } \alpha_x > 0; \mbox{  or} \nonumber\\
	& \alpha_x < \alpha_t \mbox{ and } \alpha_x < 0, \mbox{ thus } \alpha_t <0 .
\end{align*}
\label{cor1-Lemma}
\end{cor}

\begin{cor}
The power index of $x$ in (\ref{gen-selfsimilarsol1}) is positive if
\begin{align*}
	& \alpha_x > \alpha_t \mbox{ and } \alpha_x > 0  \mbox{  or} \nonumber\\
	& \alpha_x < \alpha_t \mbox{ and } \alpha_x < 0.
\end{align*}
\label{cor2-Lemma}
\end{cor}

Regarding the existence of our chosen solutions: we begin by assuming that the NSE has solutions, smooth in $({\rm I\!R}^3 \times 0)$; this is the essence of Conjecture \ref{Conjecture} and it is done to ensure the fluid has finite initial energy. We can also assume our chosen solutions are nonzero at $t=0$, because the initial velocity and pressure are the energy source of the system. 

The NSE admits the scaling transformation (\ref{EK-scaling}). Suppose all NSE solutions that scale-transform according to (\ref{EK-scaling}) are discontinued at $t=0$. Such supposition would be incorrect (e.g., see the initial solution for the Taylor-Green vortex, eq.(12) in \cite{Tay37}). We will therefore suppose, that based on Conjecture \ref{Conjecture}, the NSE always has nonzero solutions, smooth in $({\rm I\!R}^3 \times 0)$; and based on the scaling symmetry of the NSE we will also assume that some of these solutions scale-transform according to (\ref{EK-scaling}). Then, they must be in the form (\ref{gen-selfsimilarsol1}), according to Lemma \ref{Lemma}.
\begin{thm}
	Bouton's self-similar solutions (\ref{gen-selfsimilarsol}) of the NSE in 3D have global in-time regularity as long as they are energy-subcritical, smooth at $t=0$ and nonzero at $t=0$. 
\label{T2}
\end{thm}

\begin{proof}
We begin by considering all of Bouton's self-similar solutions (\ref{gen-selfsimilarsol}). To keep the NSE invariant under scaling transformation, the velocity and energy transform as
\begin{eqnarray*}
	\vec{u}^{\prime}&=&k^{\alpha_x-\alpha_t}\vec{u},  \nonumber \\
	E^{\prime}&=&k^{5\alpha_x-2\alpha_t}E.
\end{eqnarray*}
The energy $E$ at a given moment $t$, would be scale-invariant if $\vec{u}$ transformed according to the law 
\begin{equation*}
	\vec{u}^\prime=k^{-\frac{3}{2}\alpha_x}\vec{u}.
\end{equation*}
Compare the transformations $k^{\alpha_x-\alpha_t}\vec{u}$ and $k^{-\frac{3}{2}\alpha_x}\vec{u}$; we have three cases to consider:
\begin{eqnarray*}
	\alpha_x-\alpha_t &<& -\frac{3}{2} \alpha_x  \nonumber \\
	 \alpha_x-\alpha_t &>& -\frac{3}{2} \alpha_x  \nonumber \\
	 \alpha_x-\alpha_t &=& -\frac{3}{2} \alpha_x.
\end{eqnarray*}
The three cases are compared whether they increase or decrease $u$ for small scales $k \rightarrow 0$. We arrive to a 5/2 law governing the criticality of the NSE:
\begin{align}
   & \frac{\alpha_t}{\alpha_x} > \frac{5}{2}      \nonumber\\
	\mbox{subcritical, when  } &\alpha_x>\alpha_t, \alpha_x<0, \alpha_t<0, \nonumber\\
	\mbox{subcritical, when  } &\alpha_x<\alpha_t, \alpha_x>0, \alpha_t>0;
\label{5/2-1}	
\end{align}
\begin{align}
   & \frac{\alpha_t}{\alpha_x} < \frac{5}{2}      \nonumber\\
	\mbox{supercritical, when  } &\alpha_x>\alpha_t, \alpha_x<0, \alpha_t<0, \nonumber\\
	\mbox{supercritical, when  } &\alpha_x<\alpha_t, \alpha_x>0, \alpha_t>0, \nonumber\\
	\mbox{supercritical, when  } &\alpha_x>\alpha_t, \alpha_x>0, \nonumber\\
	\mbox{subcritical, when  } &\alpha_x<\alpha_t, \alpha_x<0, \nonumber\\
	\mbox{subcritical, when  } &\alpha_x=\alpha_t;
\label{5/2-2}
\end{align}
\begin{align}
   & \frac{\alpha_t}{\alpha_x} = \frac{5}{2}      \nonumber\\
	\mbox{critical, when  } &\alpha_x>\alpha_t, \alpha_x<0, \alpha_t<0, \nonumber\\
	\mbox{critical, when  } &\alpha_x<\alpha_t, \alpha_x>0, \alpha_t>0.
\label{5/2-3}
\end{align}
	The case $\alpha_x=1$, $\alpha_t=2$ yields $\alpha_t/\alpha_x=2 < 5/2$, supercritical. We omit the possible combinations that involve $\alpha_x=0$ since in them the $x,y,z$ coordinates do not rescale. Now we will use the provisions of Conjecture \ref{Conjecture} and will limit our study only to those self-similar solutions of Bouton, that are not equivalent to zero and are smooth in $({\rm I\!R}^3 \times 0)$. Such solutions will have the form (\ref{gen-selfsimilarsol1}); the corollaries of  Lemma \ref{Lemma} show the conditions that need to be met. 

\begin{prop}
	According to Corollary \ref{cor2-Lemma}, to have solutions smooth in $x$ we need $\alpha_x > \alpha_t, \alpha_x > 0$, irrespective of the sign of $\alpha_t$. This is a supercritical case by virtue of (\ref{5/2-2}). If $\alpha_t>0$, the solutions are smooth in $t$, according to Corollary \ref{cor1-Lemma}.  
\label{prop1}
\end{prop}

\begin{prop}
	According to Corollary \ref{cor2-Lemma}, we may also have solutions smooth in $x$, as long as $\alpha_x < \alpha_t, \alpha_x < 0$, irrespective of the sign of $\alpha_t$. This is a subcritical case by virtue of (\ref{5/2-2}). If $\alpha_t<0$, the solutions are smooth in $t$, according to Corollary \ref{cor1-Lemma}.  
\label{prop2}
\end{prop}

\begin{prop}
According to Corollary \ref{cor1-Lemma}, to have solutions smooth in both $t$ and $x$ we need $\alpha_x > \alpha_t, \alpha_t > 0, \alpha_x > 0$; this is a supercritical case by virtue of (\ref{5/2-2}), irrespective of the sign of $\alpha_t$.
\label{prop3}
\end{prop}

\begin{prop}
According to Corollary \ref{cor1-Lemma}, we may also have solutions smooth in both $t$ and $x$, as long as $\alpha_x < \alpha_t, \alpha_t < 0, \alpha_x < 0$; this is a subcritical case by virtue of (\ref{5/2-2}).
\label{prop4}
\end{prop}
	Consider the solutions of Proposition \ref{prop4}; they are energy-subcritical and smooth in $({\rm I\!R}^3 \times [0,\infty))$. Suppose they are diverging at some time $t \in (0, T]$. According to the Beale-Kato-Majda criterion \cite{BKM},
\begin{align*}
	& \int_0^T \mbox{sup}|\nabla \times \vec{u}| dt = \infty; \nonumber \\
	& \Rightarrow \mbox{   } \exists t: |\nabla \times \vec{u}| = \infty; \nonumber \\
	& \Rightarrow \mbox{   } \exists i: \frac{\partial \vec{u}}{\partial x_i} = \infty, \nonumber
\end{align*}
	where $i = 1,2,3$ and $x_1=x; x_2=y; x_3=z$. Suppose this is induced solely by rescaling, that is, a smooth solution of the NSE begins to diverge due to the scaling invariant properties of the equation, 
\begin{align*}
	 \left( \frac{\partial \vec{u}}{\partial x_i} \right)^\prime = \frac{k^{\alpha_x-\alpha_t}}{k^{\alpha_x}} \frac{\partial \vec{u}}{\partial x_i} = \nonumber \\
	= k^{-\alpha_t} \frac{\partial \vec{u}}{\partial x_i} = \frac{1}{k^{\alpha_t}}  \frac{\partial \vec{u}}{\partial x_i} = \infty, 
\end{align*}
where we take $(\partial \vec{u}/\partial x_i)^{\prime} = (\partial \vec{u^{\prime}}/\partial x^{\prime}_i)$, just as Bouton takes $y^{\prime}_1 = dy_1/dx_1$ in eq. (12) in \cite{Bou99} (in his work, the index 1 denotes transformation).
At fine scales, $k \rightarrow 0$ and therefore it must be true that 
\[
	\alpha_t >0, 
\]
which is a contradiction, since according to Proposition \ref{prop4}, our chosen solutions require $\alpha_t<0$. 
\end{proof}

Therefore, rescaling does not lead to the blow-up of our chosen solutions. This can be shown even for the supercritical cases. Consider the supercritical cases of Propositions \ref{prop1}, \ref{prop3}. These solutions are discontinued at $t=0$, but are smooth everywhere else. The proof in Theorem \ref{T2} applies here as well (excluding the point $t=0$), since again $\alpha_t<0$. Since our chosen solutions are isobaric polynomials in $x,y,z$ and $t$, the Beale-Kato-Majda criterion confirmed that their smoothness in $({\rm I\!R}^3 \times [0,\infty))$ cannot lead to divergence through rescaling.

\section{Bounded Energy}

Consider the energy of the fluid for the general scaling (\ref{EK-scaling}); the energy balance is given by \cite{Wie19}
\begin{eqnarray*}
	\frac{1}{2}\int _{\rm I\!R^3} |\vec{u}(x,y,z,t)|^2 dxdydz & + & \nu \int_0^t \int _{\rm I\!R^3} |\nabla \vec{u}(x,y,z,s)|^2 dxdydzds = \nonumber \\ 
	& = & \frac{1}{2}\int _{\rm I\!R^3} |\vec{u}(x,y,z,0)|^2 dxdydz.
	\label{ebalance}
\end{eqnarray*}
The total energy will remain finite through energy conservation: we are guaranteed that $\vec{u}_0 \in C^\infty ({\rm I\!R}^3)$ and expect both the Rhs and Lhs to be finite, as long as the energy is subcritical as per Theorem \ref{T2} which ensures that the velocity does not diverge anywhere in ${\rm I\!R^3}$ and for $\forall t$. 


\section{The Second Law}

The question has been posed as to whether the second law of thermodynamics may be violated by ordering the system through rearranging of its energy \cite{Tao07}. The NSE (\ref{NSE}) models an incompressible, isothermal fluid; it is therefore isentropic, while energy arrangements such as raising $E$ in one place while decreasing it in another would mean a decrease in entropy. This is forbidden by definition as the entropy of the fluid is $S=C_{\mbox{average}} \ln T = \mbox{const}$ because $T=\mbox{const}$. The possible change in entropy therefore pose a physical problem rather than a mathematical one. In attempting to define $T$ for a continuous medium we recognize that the isothermal, incompressible NSE which models continuous medium does not contain enough physics in order to guarantee $S=\mbox{const}$. Rather, it will be assumed for now that entropy does not change the same way this is assumed for the temperature. We will consider this as an {\it a priori} imposed condition on the fluid system rather than a corollary. 

In comparison, energy conservation in the fluid system is expected to hold from a physics' standpoint; it can be approached mathematically as well, since the first vectorial equation within the NSE imposes conservation of momentum, which would guarantee that energy will be conserved. However, no mathematical condition on the temperature is imposed in the NSE (\ref{NSE}) and thus $T$ cannot be controlled mathematically and it can be only supposed that the second law of thermodynamics would not be violated.

\section{When Viscosity does not scale}

In this section, we will work with the scaling transformation (\ref{scaling}), where the viscosity $\nu$ is kept at the same value at all length scales. This is the case $\alpha_x=1$, $\alpha_t=2$ (or, $\alpha_t=2\alpha_x$) in (\ref{EK-scaling}). Bouton's theorem (\ref{T2-Bouton}) yields Leray's self-similar solutions \cite{Can96}
\begin{eqnarray}
	\vec{u}(x,y,z,t) = \frac{1}{\sqrt{t}} \mathbf{F}\left( \frac{x}{\sqrt{t}},\frac{y}{\sqrt{t}},\frac{z}{\sqrt{t}} \right) \nonumber \\
	p(x,y,z,t)= \frac{1}{t} F \left( \frac{x}{\sqrt{t}},\frac{y}{\sqrt{t}},\frac{z}{\sqrt{t}} \right).
	\label{Leray}
\end{eqnarray}
In the literature, Jean Leray \cite{Ler34} is credited with noticing that due to its scaling invariance the NSE admits the so-called {\it self-similar} solutions \cite{Erc15}, \cite{Brad18} $e^{\varepsilon} \vec{u}(e^{\varepsilon}(x,y,z),e^{2\varepsilon}t)$, $e^{2\varepsilon}p(e^{\varepsilon}(x,y,z),e^{2\varepsilon}t)$, where $\varepsilon \ge 0 $ \cite{Fus03}. According to Bouton's Theorem \ref{T2-Bouton} for the NSE, Leray's self-similar solutions are the only possible self-similar solutions which the scaling transformation (\ref{scaling}) can generate. They are undefined at $t=0$ or at the origin etc., resulting in divergence of the fluid energy at the initial moment or at the origin etc. (contrary to the provisions of Conjecture \ref{Conjecture}). Analogous reasoning applies to the general scaling transformation (\ref{EK-scaling}) with $\alpha_t=2\alpha_x$, where again the solutions have fluid energy undefined at $t=0$ or at the origin/along the axes.

Continuing our study of the case of non-scaling viscosity $\alpha_t=2\alpha_x$, we note that Bouton's theorem (\ref{T2-Bouton}) can be changed to take units, or dimensionalization, into account:
\begin{thm} (of Bouton): In the case $\alpha_t=2\alpha_x$, in order to meet all isobarity requirements of their scaling transformation, the relative invariant functions $u(x,y,z,t)$, $v(x,y,z,t)$ and $w(x,y,z,t)$ must be homogeneous in $x,y,z$ of degree 1 and homogeneous in $t$ of degree -1; while the relative invariant function $p(x,y,z,t)$  must be homogeneous in $x,y,z$ of degree 2 and homogeneous in $t$ of degree -2. The Navier-Stokes system of PDEs becomes: 
\begin{eqnarray}
	\frac{\partial \vec{u}}{\partial t} + (\vec{u} \cdot \nabla)\vec{u} &=& \nu  \Delta \vec{u} - \nabla p \nonumber\\
	\nabla \cdot \vec{u} &=& 0  \nonumber\\
	(\vec{r}\cdot\nabla)\vec{u}&=&\vec{u}  \nonumber\\
	 t\frac{\partial \vec{u}}{\partial t}&=& -  \vec{u} \nonumber\\
	(\vec{r}\cdot\nabla)p&=&2p \nonumber\\
	t\frac{\partial p}{\partial t} & = & -2p .
	\label{NSE2}	
\end{eqnarray}
\label{T3-Bouton}
\end{thm}

The added equations represent Bouton's conditions of homogeneity. They yield a set of solutions, which are a subset of Leray's self-similar solutions (\ref{Leray}):
\begin{eqnarray}
	\vec{u}(x,y,z,t) &=& \frac{x}{t} \mathbf{F}\left(\frac{y}{x},\frac{z}{x}\right)  \nonumber \\
	p(x,y,z,t)&=& \frac{x^2}{t^2} F\left(\frac{y}{x},\frac{z}{x}\right),
	\label{dim-sol}
\end{eqnarray}
where again the solutions have fluid energy undefined at $t=0$. The solutions (\ref{Leray}) and (\ref{dim-sol}) can be considered defined for $t>0$ and $t<0$ \cite{Brad18}. We note that the self-similar solutions (\ref{dim-sol}) of the PDE system (\ref{NSE2}) are in the form $\varphi(t)F(x,y,z)$ - the time dependence is a factor which multiplies a function of the spatial coordinates. (\ref{dim-sol}) can be considered to have proper units, if this is necessary for a study of the dimensionalized NSE. Through homogeneity, $u(x,y,z,t)$, $v(x,y,z,t)$, $w(x,y,z,t)$ are isobaric IRFs of weight $W=-1$; while $p(x,y,z,t)$ is an isobaric IRF of weight $W=-2$. Proper isobarity in both $x,y,z$ and $t$ is a necessary and sufficient condition for their correct transformation.

\begin{thm}
	Within the set of Leray's self-similar solutions (\ref{Leray}), there is always a subset (\ref{dim-sol}) which is not subject to scaling-induced blow-up.
\label{T3}
\end{thm}
\begin{proof}
According to \cite{Llo81}, the NSE admits the following Lie groups of transformations:

Time translations
\begin{eqnarray}
	\mathcal{T}&=&\frac{\partial}{\partial t} \nonumber \\
	t' &=& t + \tau;  \nonumber \\
	-\infty < \tau < \infty,& &  \tau=\mbox{const}
    \label{timetranslations}
\end{eqnarray}

Scale changes
\begin{eqnarray}
	\mathcal{X} = x \frac{\partial}{\partial x} + y \frac{\partial }{\partial y} + z \frac{\partial }{\partial z} + 2t \frac{\partial }{\partial t} 
       - u \frac{\partial }{\partial u} - v \frac{\partial }{\partial v} - w \frac{\partial }{\partial w} - 2p \frac{\partial }{\partial p}
	\label{scaling-oper}	
\end{eqnarray}

Infinitesimal rotations

\begin{eqnarray}
	\mathcal{R}_x= y\frac{\partial}{\partial z}-z\frac{\partial}{\partial y}+v\frac{\partial}{\partial w}-w\frac{\partial}{\partial v}; \mbox{  } 
	\mathcal{R}_y= z\frac{\partial}{\partial x}-x\frac{\partial}{\partial z}+w\frac{\partial}{\partial u}-u\frac{\partial}{\partial w}; \nonumber\\
	\mathcal{R}_z= x\frac{\partial}{\partial y}-y\frac{\partial}{\partial x}+u\frac{\partial}{\partial v}-v\frac{\partial}{\partial u}. 
    \label{rotations}
\end{eqnarray}

Because the operators (\ref{timetranslations})-(\ref{rotations}) do not depend on partial derivatives, all absolute covariants of the NSE will depend only on $u,v,w,p$. This is a property of the NSE worth noticing; its importance will be seen in the last chapter of this work.

The necessary and sufficient condition (\ref{XI}) the absolute invariants/covariants of the NSE must satisfy is expressed in the system of five PDEs
\begin{eqnarray}
      \mathcal{T}I = 0 \nonumber \\
      \mathcal{X}I = 0 \nonumber \\
      \mathcal{R}_x I = 0 \nonumber \\
      \mathcal{R}_y I = 0 \nonumber \\
      \mathcal{R}_z I = 0. 
      \label{0-form-system}
\end{eqnarray}
	We solve (\ref{0-form-system}) with Maple \cite{Map19} and obtain all absolute covariants
\begin{equation}
	I = F \left( \frac{p}{u^2+v^2+w^2} \right),
	\label{allcovariants}
\end{equation}
	where $F$ is an arbitrary function and $p/\rho$ is written shortly as ``$p$''. It needs to be noted, that the form of all covariants (\ref{allcovariants}) is valid under the transformation (\ref{EK-scaling}) as well as under the general scaling (\ref{EK-scaling}); in the first case we solve the PDE system (\ref{0-form-system}), while in the later $\mathcal{X}$ is replaced by  two scaling operators. In the simplest case (\ref{allcovariants}) equals the dimensionless Euler's number 
\begin{equation}
	\mathcal{E} = \frac{p}{u^2+v^2+w^2},
	\label{Eu}
\end{equation}
also known as cavitation number (if the vapor pressure is set to zero). $\mathcal{E}$ is an invariant with respect to all symmetry transformations of the NSE (eqs. (\ref{0-form-system})). From (\ref{dim-sol}) it follows that due to the factorization of the self-similar solution, the time-dependent parts in $\mathcal{E}$ cancel out and it turns out to be a conserved quantity (in time),
\begin{eqnarray}
	\mathcal{E} = \frac{\frac{x^2}{t^2} F_4(\frac{y}{x},\frac{z}{x})}{\frac{x^2}{t^2}\Bigl[ F_1^2(\frac{y}{x},\frac{z}{x})+F_2^2(\frac{y}{x},\frac{z}{x})+F_3^2(\frac{y}{x},\frac{z}{x})\Bigr] }= \frac{F_4}{F_1^2+F_2^2+F_3^2}=\mbox{const}.
	\label{Eu-const}
\end{eqnarray}
	If the conserved quantities are found to provide upper bounds to the size of the solution, they are termed {\it coercive} \cite{Tao07}. The cavitation number $\mathcal{E}$ contains the square of the modulus of the velocity $|\vec{u}|^2$ in denominator and can thus be considered coercive; it is not an integral over ${\rm I\!R^3}$ and thus not subject to delta-functions within the integration domain. Suppose that at large scales $k\gg1$ the modulus of the velocity is finite but the energy of the system cannot ensure it does not tend to infinity at a certain location at very fine scales. $\mathcal{E}$ is scale-invariant which means that as long as it has finite value at large scales it must retain them at arbitrarily fine scales as well. This would guarantee that $|\vec{u}|^2$ does not grow uncontrollably at any locale of infinitely small size. In addition, this control is being exercised at any given moment in time, because $\mathcal{E}$ is conserved, (\ref{Eu-const}). Also, $|\vec{u}|$ and the pressure $p$ cannot be diverging together as per the Bernoulli principle for incompressible fluids; and therefore, the ratio $p/|\vec{u}^2|$ cannot be finite unless both its numerator and denominator are finite $\forall t$.	
\end{proof}
\begin{thm}
	When the viscosity is not allowed to scale, $\alpha_x=1$, $\alpha_t=2$ (or, $\alpha_t=2\alpha_x$), the Navier-Stokes equation has a set of non self-similar solutions, which are nonzero at $t=0$, nontrivial (not identically zero everywhere) and infinitely differentiable throughout ${\rm I\!R}^3$ and $\forall t$.
\end{thm}
\begin{proof}
	According to Lemma \ref{Lemma}, The NSE exhibits self-similar solutions (\ref{gen-selfsimilarsol1}) under the general scaling transformation (\ref{EK-scaling}) when viscosity is allowed to scale. If the viscosity is not allowed to scale i. e.,  $\alpha_x=1$, $\alpha_t=2$ (or, $\alpha_t=2\alpha_x$), the solutions (\ref{gen-selfsimilarsol1})
\begin{eqnarray*}
	\vec{u} &=&\mathbf{C} t^{\frac{\beta_x-\beta_t}{\beta_t}} +  x^{\frac{\beta_x-\beta_t}{\beta_x}} \mathbf{P} \left( \frac{y}{x},\frac{z}{x}  \right)  \nonumber \\
	p &=& C t^{\frac{2(\beta_x-\beta_t)}{\beta_t}}  +  x^{\frac{2(\beta_x-\beta_t)}{\beta_x}} P \left( \frac{y}{x},\frac{z}{x}  \right), \nonumber \\
	 \beta_t &\ne & 2 \beta_x
\end{eqnarray*}
	still exist (as per the discussion following Corollary \ref{cor2-Lemma}), however they do not transform according to 
\begin{eqnarray*}	
	(u',v',w')&=&(1/k)(u,v,w) \nonumber \\
	p'&=&(1/k^2)p,  \nonumber \\
	0<k<\infty,& &k=\mbox{const}, 
\end{eqnarray*}
	and are therefore non self-similar under this transformation. Consider the supercritical (under the general scaling transformation (\ref{EK-scaling}))  solutions of Proposition \ref{prop3}: here, $\beta_x>\beta_t$ and $\beta_x > 0$. We will set $\beta_t>0$, since its sign is irrelevant. Suppose that the solutions are smooth and then diverging after finite time. Suppose this is induced solely by rescaling, that is, a smooth solution of the NSE begins to diverge due to the scaling invariant properties of the equation. According to the Beale-Kato-Majda criterion, the integral of the vorticity will diverge and  since $\alpha_t=2\alpha_x$, it must be true that
\begin{equation*}
	\left( \frac{\partial \vec{u}}{\partial x_i} \right)^\prime = k^{\beta_x-\beta_t}\frac{\left[k^{\frac{2}{\beta_t}}\mathbf{C}t+k^{\frac{1}{\beta_x}}x\mathbf{P}(\frac{y}{x},\frac{z}{x}) \right]}{\partial x_i} = \infty 
\end{equation*}
	at fine scales $k \rightarrow 0$, which would be possible as long as $\beta_x<\beta_t$, or $\beta_x < 0$ or $\beta_t<0$. However, any of these would contradict the supercriticality condition of Proposition \ref{prop3}. It is not possible therefore that this set of solutions would exhibit scaling induced blow-up in finite time. Moreover, since the non self-similar solutions of Bouton are nonzero at $t=0$ and nontrivial (not identically zero everywhere) as per Lemma \ref{Lemma}, by virtue of the same Lemma \ref{Lemma} they are polynomials, and therefore infinitely differentiable in ${\rm I\!R}^3 \times [0, \infty)$.   
\end{proof}

\section{Turbulence}
We calculate the scaling invariants of the NSE in the case $\alpha_x=1, \alpha_t=2$ (transformation (\ref{scaling})), independently of other transformations by writing $\mathcal{X}I=0$, where the infinitesimal operator of scale changes $\mathcal{X}$ is given in (\ref{scaling-oper}). Maple yields 
\begin{equation*}
	I = F \left(\frac{y}{x},\frac{z}{x},ux,vx,wx,px^2,\frac{t}{x^2} \right), 
\end{equation*}
$F$ being an arbitrary function. Then $|\vec{u}|^2 r^2$ is an absolute covariant ($r = |\vec{r}\mbox{ }| = \sqrt{x^2+y^2+z^2}$); so are $(\vec{u}\cdot\vec{r})^2$, $(\vec{r} \times \vec{u})_i$ and $(\vec{r} \times \vec{u})^2$ etc. Write $I$ as
\begin{equation}
	|\vec{r} \times \vec{u}| + \left(\frac{1}{\nu}\right) pr^2.
	\label{ang-mom-inv}
\end{equation}
This is a scale-invariant quantity in which $|\vec{r} \times \vec{u}|$ is the modulus of the angular momentum per unit mass for fluid at position $\vec{r}$; the pressure-viscosity term $(1/\nu) pr^2$ has the same units and can be considered to account for viscous transfer of angular momentum. Next, consider the total angular momentum of the fluid per unit mass:
\begin{equation}
	\frac{\vec{H}_{\mbox{tot}}}{m_0} =\sum\limits_{\rm I\!R^3}  (\vec{r} \times \vec{u}) =\mbox{const}.
	\label{Htot}
\end{equation}
The RHS is scale-invariant vector because $X(\vec{r} \times \vec{u})_i=0$; then the total angular momentum per unit mass is a scale-invariant and a conserved quantity as well (recall that $m_0=\mbox{const}$ is the unit mass of the fluid particle at position $\vec{r}$ as fluid density is the same everywhere). 

It may be argued, that the above mentioned scale-invariants as well as eqs. (\ref{ang-mom-inv}) and (\ref{Htot}) suggest collective fluid behavior: the scaling invariance of $\vec{r} \times \vec{u}$ and $\vec{r} \cdot \vec{u}$ means that if at some scale $|\vec{u}|$ decreases yet at a different scale $|\vec{r} \times \vec{u}|$ will remain invariant if the angle between $\vec{r}$ and $\vec{u}$ changes ($|\vec{r} \times \vec{u}|=|\vec{r}||\vec{u}|\sin \gamma$). The invariant behavior of $\vec{r} \cdot \vec{u}$ suggests the same: the change of $|\vec{u}|$ at one scale may trigger changes of the direction of fluid motion at other scales ($\vec{r} \cdot \vec{u}=|\vec{r}||\vec{u}|\cos \gamma$). Small changes in $|\vec{u}|$ can set off large changes in the directions of motion at smaller scales. Nonetheless, such momentum exchanges between scales can only take place in different moments because time is scaled as well. 

However, (\ref{Htot}) represents a conserved quantity. Through scale invariance, any changes in $\vec{r} \times \vec{u}$ on a large scale will trigger instantaneous changes in the small scales, where rotation can take place. This is a collective behavior throughout the fluid volume, made possible by viscous transfer of momentum (i.e. the second term in (\ref{ang-mom-inv})) and ultimately due to the scaling invariance of the NSE.

Then, turbulence can generally be looked upon as the collective tendency of a viscous fluid to generate small scale vortices via its scaling properties and exchange of momentum. Such explanation is only qualitative though, while what is desired is a quantitative description which would allow control of turbulent phenomena.

\section{Differential Forms as Conserved Quantities}

Bouton's work consists of deriving the invariants of a differential equation; this is accomplished by examining the finite transformations admitted by the equation as well as by putting to use their infinitesimal operators to find absolute and relative invariants. Both approaches give the same results; the infinitesimal transformations in equations (30) and (40) in Bouton's article \cite{Bou99} are published by him for the first time; they confirm his first Theorem, which is based on finite transformations. In the sections above, we applied Bouton's analysis to the NSE. In this section, we will limit our study to the scaling transformation (\ref{scaling}), where $\alpha_x=1,\alpha_t=2$, or $\alpha_t=2\alpha_x$, since it yields a conserved quanity, $\mathcal{E}$. Recall equation (\ref{0-form-system}): all transformations, admitted by the NSE were used together to form the system of 5 PDEs $\mathcal{T}I=0$, $\mathcal{X}I=0$,... in order to find a common covariant. As a matter of fact, the PDE system (\ref{0-form-system}) is equivalent to the system of Lie derivatives
\begin{eqnarray}
   \mathcal{L}_{\mathcal{T}}\omega^n=0 \nonumber\\
   \mathcal{L}_{\mathcal{X}}\omega^n=0 \nonumber\\
   \mathcal{L}_{\mathcal{R}_x}\omega^n=0 \nonumber\\
   \mathcal{L}_{\mathcal{R}_y}\omega^n=0 \nonumber\\
   \mathcal{L}_{\mathcal{R}_z}\omega^n=0,
   \label{n-form}
\end{eqnarray}
where $I = \omega^0$ of eq.(\ref{0-form-system}) is a zero-form covariant (when $n=0$) along the flows of the vector fields $\mathcal{T},\mathcal{X},\mathcal{R}_x,\mathcal{R}_y,\mathcal{R}_z$. It makes sense to study precisely these flows, since they correspond to vector fields that leave the NSE unchanged (i.e., if we move all terms of the NSE to the Lhs and call this $L$, it is said that the NSE admits the transformation if each operator, applied to $L$, reduces it to zero whenever $L=0$. For an example application of this rule see \cite{Llo81}). Therefore, Bouton's invariant analysis of the NSE is but a first stage from a total of nine stages, all of which yield invariants of the NSE. The nine stages correspond to the 8 variables $x,y,z,t,u,v,w,p$ with whose differentials we can build 0-, 1-, ... and 8-forms. In all stages we use Maple to calculate the solutions of the corresponding systems of PDEs.

In all cases, we calculate (\ref{n-form}) seeking $n-$form covariants as expressions of the dependent variables $u,v,w,p$. Then we proceed as follows: we replace these variables with their corresponding functions, the solutions of the NSE $u(x,y,z,t),v(x,y,z,t),w(x,y,z,t)$ and $p(x,y,z,t)$. We then examine the resulting expressions and focus on those of them, which are conserved in time due to the cancelling out of the time-dependent factor. From Theorem \ref{T3-Bouton} it follows that they will also be scale-invariant. The fractions $p/\vec{u}^2$, $u/v$, $|\vec{u}|/v$, $\vec{u}^2/(uv)$ etc. are time-independent by virtue of (\ref{dim-sol}). This result follows from Theorem \ref{T3-Bouton}, according to which the self-similar solutions contain time as a factor and thus the above ratios represent dimensionless scalar fields that no longer depend on $t$. In all expressions, $\vec{u}^2=u^2+v^2+w^2$; $|\vec{u}|=\sqrt{u^2+v^2+w^2}$. 

The idea to utilize Cartan's calculus in the study of the NSE has been researched previously by R. Kiehn \cite{Kie03}, \cite{Kie01}. In his analysis, he has studied the Lie derivatives along the flows of vector fields and their application to the NSE while in this work, the focus is on the vector fields of the differential operators admitted by the NSE. The significance of the derived results below is subject to future research where it would be determined how the invariants of the NSE are related to the physics of the flow.

{\underline {\it 0-forms:}} To find the 0-forms that are covariants of all five operators simultaneously, we write eq. (\ref{0-form-system}). As explained, this is equivalent to the system of 5 Lie derivatives (\ref{n-form}) when $n=0$; it is a system of 5 PDEs, which yields the covariants (\ref{allcovariants}). In the simplest case, this is (\ref{Eu}), the Euler's, or, cavitation number $\mathcal{E}$. Thus, $B_0=\mathcal{E}$ is an absolute invariant 0-form of the NSE; time-independent by virtue of (\ref{dim-sol}) and thus a conserved quantity of Bouton $B_0$. It is a dimensionless scalar field in ${\rm I\!R}^3$; its name ``cavitation number'' indicates that $\mathcal{E}$ is related to cavitation phenomena and thus to turbulence. All conserved quiantities below are multiples of $\mathcal{E}$.

Recall again the fact that all infinitesimal operators in (\ref{timetranslations})-(\ref{rotations}) do not depend on the derivatives of $u,v,w,p$ as demonstrated by Lloyd \cite{Llo81}. This is an advantageous circumstance based on which all NSE invariants are expressions of the dependent variables alone. Because they do not contain derivatives, we can use them to find conserved quantities in a simple functional form.  

{\underline {\it 1-forms:}} Let $\omega^1$ be the form $\omega^1=a_1(u,v,w,p)dx+a_2(u,v,w,p)dy+...+a_8(u,v,w,p)dp$. We solve the system (\ref{n-form}) for $k=1$, which contains $8\times5 =40$ PDEs. No conserved quantities are found. All functions $a_i$ are considered to contain the dependent variables $u,v,w,p$. 

{\underline {\it 2-forms:}} Let $\omega^2$ be the form $\omega^2=a_1(u,v,w,p)dx \wedge dy+a_2(u,v,w,p)dx \wedge dz + ...+a_{28}(u,v,w,p)dw \wedge dp$. The form can be reduced to 28 terms. We solve the system (\ref{n-form}) for $k=2$, which contains $28\times5 =140$ PDEs. We obtain a 2-form conserved quantity of Bouton $B_2$, which has 10 terms. All are multiplied by $F(B_0)$, where $F$ is an arbitrary function. For simplicity, we set $F(B_0) = B_0=\mathcal{E}$ and get

\begin{eqnarray*}
	B_2 = \mathcal{E} \Bigl[ \left(\frac{\vec{u}^2+u^2}{\vec{u}^2} \right) dx \wedge du + \Bigr. 
	                  \left(\frac{\vec{u}^2w+u|\vec{u}|v}{\vec{u}^2|\vec{u}|} \right) dx \wedge dv + 
			  \left(\frac{u|\vec{u}|w-v\vec{u}^2}{\vec{u}^2|\vec{u}|} \right) dx \wedge dw + \\
                  	+ \left(\frac{u|\vec{u}|v-w\vec{u}^2}{\vec{u}^2|\vec{u}|} \right) dy \wedge du +
	                  \left(\frac{\vec{u}^2+v^2}{\vec{u}^2} \right) dy \wedge dv + 
			  \left(\frac{|\vec{u}|v^3+u\vec{u}^2w}{\vec{u}^2|\vec{u}|w} \right) dy \wedge dw + \\
	                + \left(\frac{u|\vec{u}|w-v\vec{u}^2}{\vec{u}^2|\vec{u}|} \right) dz \wedge du + 
			  \left(\frac{|\vec{u}|vw-u\vec{u}^2}{\vec{u}^2|\vec{u}|} \right) dz \wedge dv + 
			  \left(\frac{\vec{u}^2+w^2}{\vec{u}^2} \right) dz \wedge dw + \\
			+ \Bigl. dt \wedge dp \mbox{ } \Bigr].
\end{eqnarray*}

{\underline {\it 3-forms:}} Let $\omega^3$ be the form $\omega^3=a_1(u,v,w,p)dx \wedge dy \wedge dz +a_2(u,v,w,p) dx \wedge dy \wedge dt +...+a_{56}(u,v,w,p)du \wedge dw \wedge dp$. The form can be reduced to 56 terms. We solve the system (\ref{n-form}) for $k=3$, which contains $56\times5 =280$ PDEs. We obtain a 3-form conserved quantity of Bouton $B_3$, which has 6 terms. All are multiplied by $F(B_0)$, where $F$ is an arbitrary function. For simplicity, we set $F(B_0) = B_0=\mathcal{E}$ and get
\begin{eqnarray*}
	B_3 = \mathcal{E} \Bigl[ \frac{w}{|\vec{u}|} dx \wedge dy \wedge dp - \frac{v}{|\vec{u}|} dx \wedge dz \wedge dp + \Bigr. 
	                  \frac{u}{|\vec{u}|} dy \wedge dz \wedge dp + \frac{w}{|\vec{u}|} dt \wedge du \wedge dv - \\
			  \Bigl. -\frac{v}{|\vec{u}|} dt \wedge du \wedge dw + \frac{u}{|\vec{u}|} dt \wedge dv \wedge dw \Bigr]. 
\end{eqnarray*}

{\underline {\it 4-forms:}} Let $\omega^4$ be the form $\omega^4=a_1(u,v,w,p)dx \wedge dy \wedge dz \wedge dt +a_2(u,v,w,p)dx \wedge dy \wedge dz \wedge du + ...+a_{70}(u,v,w,p)du \wedge dv \wedge dw \wedge dp$. The form can be reduced to 70 terms. We solve the system (\ref{n-form}) for $k=4$, which contains $70\times5 =350$ PDEs. We obtain a 4-form conserved quantity of Bouton $B_4$, which has 18 terms. All are multiplied by $F(B_0)$, where $F$ is an arbitrary function. For simplicity, we set $F(B_0) = B_0=\mathcal{E}$ and get
\begin{eqnarray*}
	B_4 = \mathcal{E} \Bigl[ \left(\frac{2\vec{u}^2-u^2-v^2}{2\vec{u}^2} \right) dx \wedge dy \wedge du \wedge dv + \Bigr.
	                  \left(\frac{2\vec{u}^2u-wv|\vec{u}|}{2\vec{u}^2|\vec{u}|} \right) dx \wedge dy \wedge du \wedge dw +\\
	                + \left(\frac{\vec{u}^2|\vec{u}|w+2uv\vec{u}^2}{2\vec{u}^2|\vec{u}|u} \right) dx \wedge dy \wedge dv \wedge dw -
	                  \left(\frac{wv|\vec{u}|+2u\vec{u}^2}{2\vec{u}^2|\vec{u}|} \right) dx \wedge dz \wedge du \wedge dv + \\
	                + \left(\frac{2\vec{u}^2-u^2-w^2}{2\vec{u}^2} \right) dx \wedge dz \wedge du \wedge dw + 
	                  \left(\frac{2\vec{u}^2uw-u^2v|\vec{u}|}{2\vec{u}^2|\vec{u}|u} \right) dx \wedge dz \wedge dv \wedge dw +\\
	                + \left(\frac{\vec{u}^2w|\vec{u}|-2uv\vec{u}^2}{2\vec{u}^2|\vec{u}|u} \right) dy \wedge dz \wedge du \wedge dv -
	                  \left(\frac{u^2v|\vec{u}|+2uw\vec{u}^2}{2\vec{u}^2|\vec{u}|u} \right) dy \wedge dz \wedge du \wedge dw +\\
                        + \left(\frac{2\vec{u}^2-v^2-w^2}{2\vec{u}^2} \right) dy \wedge dz \wedge dv \wedge dw +
                          \left(\frac{\vec{u}^2+u^2}{\vec{u}^2} \right) dx \wedge dt \wedge du \wedge dp +\\
                        + \left(\frac{w\vec{u}^2+uv|\vec{u}|}{\vec{u}^2|\vec{u}|} \right) dx \wedge dt \wedge dv \wedge dp +
                          \left(\frac{u|\vec{u}|w^2-\vec{u}^2vw}{\vec{u}^2|\vec{u}|w} \right) dx \wedge dt \wedge dw \wedge dp +\\
                        + \left(\frac{u|\vec{u}|v-\vec{u}^2w}{\vec{u}^2|\vec{u}|} \right) dy \wedge dt \wedge du \wedge dp + 
                          \left(\frac{\vec{u}^2+v^2}{\vec{u}^2} \right) dy \wedge dt \wedge dv \wedge dp +\\
                        + \left(\frac{|\vec{u}|vw^2+\vec{u}^2uw}{\vec{u}^2|\vec{u}|w} \right) dy \wedge dt \wedge dw \wedge dp +
	                  \left(\frac{u|\vec{u}|w^2+\vec{u}^2vw}{\vec{u}^2|\vec{u}|w} \right) dz \wedge dt \wedge du \wedge dp +\\
	                + \left(\frac{|\vec{u}|vw-\vec{u}^2u}{\vec{u}^2|\vec{u}|} \right) dz \wedge dt \wedge dv \wedge dp +
			  \Bigl. \left(\frac{\vec{u}^2+w^2}{\vec{u}^2} \right) dz \wedge dt \wedge dw \wedge dp \mbox{ } \Bigr].
\end{eqnarray*}

{\underline {\it 5-forms:}} Let $\omega^5$ be the form $\omega^5=a_1(u,v,w,p)dx \wedge dy \wedge dz \wedge dt \wedge du +a_2(u,v,w,p)dx \wedge dy \wedge dz \wedge dt \wedge dw + ...+a_{56}(u,v,w,p)dt \wedge du \wedge dv \wedge dw \wedge dp$. The form can be reduced to 56 terms. We solve the system (\ref{n-form}) for $k=5$, which contains $56\times5 =280$ PDEs. We obtain a 5-form conserved quantity of Bouton $B_5$, which has 6 terms. All are multiplied by $F(B_0)$, where $F$ is an arbitrary function. For simplicity, we set $F(B_0) = B_0=\mathcal{E}$ and get
\begin{eqnarray*}
	B_5=\mathcal{E} \Bigl[ \frac{u}{|\vec{u}|} dx \wedge dy \wedge dz \wedge du \wedge dp + \Bigr.
	                  \frac{v}{|\vec{u}|} dx \wedge dy \wedge dz \wedge dv \wedge dp + 
	                  \frac{w}{|\vec{u}|} dx \wedge dy \wedge dz \wedge dw \wedge dp + \\
	                + \frac{u}{|\vec{u}|} dx \wedge dt \wedge du \wedge dv \wedge dw + 
	                  \frac{v}{|\vec{u}|} dy \wedge dt \wedge du \wedge dv \wedge dw + 
			  \Bigl. \frac{w}{|\vec{u}|} dz \wedge dt \wedge du \wedge dv \wedge dw \Bigr]. 
\end{eqnarray*}
{\underline {\it 6-forms:}} Let $\omega^6$ be the form $\omega^6=a_1(u,v,w,p)dx \wedge dy \wedge dz \wedge dt \wedge du \wedge dv +a_2(u,v,w,p)dx \wedge dy \wedge dz \wedge du \wedge dt \wedge dw + ...+a_{28}(u,v,w,p)dy \wedge dt \wedge du \wedge dv \wedge dw \wedge dp$. The form can be reduced to 28 terms. We solve the system (\ref{n-form}) for $k=6$, which contains $28\times5 =140$ PDEs. We obtain a 6-form conserved quantity of Bouton $B_6$, which has 10 terms. All are multiplied by $F(B_0)$, where $F$ is an arbitrary function. For simplicity, we set $F(B_0) = B_0=\mathcal{E}$ and get
\begin{eqnarray*}
	B_6 = \mathcal{E} \Bigl[ dx \wedge dy \wedge dz \wedge du \wedge dv \wedge dw + \Bigr.
                          \left(\frac{2\vec{u}^2-u^2-v^2}{2\vec{u}^2} \right) dx \wedge dy \wedge dt \wedge du \wedge dv \wedge dp + \\
	                + \left(\frac{2\vec{u}^2u-|\vec{u}|vw}{2\vec{u}^2|\vec{u}|} \right) dx \wedge dy \wedge dt \wedge du \wedge dw \wedge dp + 
	                  \left(\frac{u^2|\vec{u}|w+2\vec{u}^2uv}{2\vec{u}^2|\vec{u}|u} \right) dx \wedge dy \wedge dt \wedge dv \wedge dw \wedge dp -\\ 
	                - \left(\frac{|\vec{u}|vw+2\vec{u}^2u}{2\vec{u}^2|\vec{u}|} \right) dx \wedge dz \wedge dt \wedge du \wedge dv \wedge dp +
                          \left(\frac{2\vec{u}^2-u^2-w^2}{2\vec{u}^2} \right) dx \wedge dz \wedge dt \wedge du \wedge dw \wedge dp - \\
	                - \left(\frac{u|\vec{u}|v-2\vec{u}^2w}{2\vec{u}^2|\vec{u}|} \right) dx \wedge dz \wedge dt \wedge dv \wedge dw \wedge dp + 
	                  \left(\frac{u|\vec{u}|w-2\vec{u}^2v}{2\vec{u}^2|\vec{u}|} \right) dy \wedge dz \wedge dt \wedge du \wedge dv \wedge dp - \\
	                - \left(\frac{u|\vec{u}|v+2\vec{u}^2w}{2\vec{u}^2|\vec{u}|} \right) dy \wedge dz \wedge dt \wedge du \wedge dw \wedge dp + 
			  \Bigl. \left(\frac{2\vec{u}^2-v^2-w^2}{2\vec{u}^2} \right) dy \wedge dz \wedge dt \wedge dv \wedge dw \wedge dp \Bigr].
\end{eqnarray*}

{\underline {\it 7-forms:}} Let $\omega^7$ be the form $\omega^7=a_1(u,v,w,p)dx \wedge dy \wedge dz \wedge dt \wedge du \wedge dv \wedge dw + ...$. We solve the system (\ref{n-form}) for $k=7$, which contains $8\times5 =40$ PDEs. No conserved quantities are found. 

{\underline {\it 8-forms:}} Let $\omega^8$ be the form $\omega^8=a_1(u,v,w,p)dx \wedge dy \wedge dz \wedge dt \wedge du \wedge dv \wedge dw \wedge dp$. Here, we have only one term. We solve the system (\ref{n-form}) for $k=8$, which contains 5 PDEs and obtain the 8-form conserved quantity of Bouton $B_8$
\begin{equation*}
   B_8=\mathcal{E} dx \wedge dy \wedge dz \wedge dt \wedge du \wedge dv \wedge dw \wedge dp.
\end{equation*}

\section{Conclusions}

The NSE global regularity conjecture is the assumption that given smooth initial conditions, a solution of the NSE can always be found such, that it will remain smooth throughout 3-dimensional space and all time. To evaluate the conjecture, we applied the pioneering work of Charles L. Bouton. According to Wikipedia, the online encyclopedia, Bouton was an American mathematician, born and lived just after the end of the American Civil War \cite{Wik22}. He obtained his Ph. D. in mathematics under the guidance of Sophus Lie at the University of Leipzig in Saxony, Germany in 1898 \cite{Osg22} and became professor in mathematics at Harvard University, where he also served as editor of leading American journals in mathematics. Bouton's invariant theory was published in 1899 \cite{Bou99}; it consists of four theorems, the most relevant to the NSE being the first one according to which the scaling invariance of the PDE requires the solutions to be isobaric functions. Bouton's work on invariants is foundational, as acknowledged by researchers \cite{Soh16}; it includes his doctoral thesis \cite{Bou99} and a number of other articles he authored. However, the memory of Bouton's life is enshrowded in mystery: despite of his prominent position in society \cite{Osg22}, which is well-documented in the media of his day, his memorial site is completely unmarked, unlike the sites pertaining to his parents and all his siblings \cite{Cot16}. His Wikipedia entry lists only gaming theory in his publications record \cite{Wik22}, while the main thrust of his work was on Lie invariants of differential equations, which remains undiscovered with very limited number of citations according to the Google Scholar database.

We adapted Bouton's first theorem to address the conjecture. Bouton's theory together with the general scaling transformation of the NSE show the general form of all self-similar solutions. Under this transformation, the NSE problem is no longer supercritical, but it can also be critical or subcritical, depending on the scaling parameters. The criticality limit is governed by a 5/2 law, determined by the ratio of two of the scaling transformation parameters. A subset of Bouton's self-similar solutions are shown to have global in-time regularity by using the Beale-Kato-Majda argument. It is assumed that such solutions always exist due to the scaling properties of the NSE and the basic physics of the flow. 

We also studied Leray's self-similar solutions as well as the regularity of a subset of those: when the solutions are homogeneous functions, they yield a number of scale-invariant conserved quantities, which are differential forms of various order and likely related to the physics of turbulence. When viscosity is kept the same at all scales, the NSE always has non self-similar, isobaric polynomial solutions smooth in ${\rm I\!R}^3 \times [0, \infty)$.

\vspace{0.5cm}
\end{document}